\documentclass{bull-l}

\usepackage{amssymb}

\begin{document}

\title[Distribution modulo one and Diophantine approximation]{ \ \ }

{  } \ \

\bigskip
\bigskip

\noindent
{\bf Distribution modulo one and Diophantine approximation, by Yann \\
Bugeaud, Cambridge University Press 2012. ISBN 978-0521111690, \\
316 pp. $\$$85}

\address{Jean-Paul Allouche \\
CNRS, Institut de Math\'ematiques de Jussieu \\
\'Equipe Combinatoire et Optimisation \\
Universit\'e Pierre et Marie Curie, Case 247 \\
4 Place Jussieu \\
F-75252 Paris Cedex 05 France}
\email{allouche@math.jussieu.fr}


 \maketitle

Let us start with the well-known fact that the base $b$ expansion ($b \geq 2$ 
an integer) of a rational number is periodic from some index on, and that the 
sequence of partial quotients of a quadratic number is periodic from some index 
on as well. These properties characterize rationals, respectively quadratic 
irrationals. On the other hand, constructing a real number with more complicated
patterns, like 
$$
0.1100010000000000000000010... = \sum_{n \geq 1} 10^{-n!}
$$ 
yields a transcendental number (Liouville number). This is also the case for 
the Barbier-Champernowne number, whose expansion in base $10$ is obtained by 
concatenating after the decimal point all the base $10$ expansions of the 
consecutive integers
$$
0.1234567891011121314151617181920...
$$
This leads to the  natural question whether the arithmetical nature of a real 
number can be inferred from its decimal expansion and vice-versa. In particular,
what can be said about the digits of $\sqrt{2}$ or about the digits of $\pi$?
Most of the elementary questions one can imagine of about the decimal digits of 
these two numbers are still unanswered, e.g., does a fixed digit (say $7$) occur
infinitely often in any of these two expansions? The nice book of Bugeaud looks 
at a whole bunch of similar questions and leads the reader to the state of the 
art, including very recent results.

\bigskip

First note that a random real number (i.e., ``almost any'' for the Lebesgue 
measure) has a quite regular distribution of its digits. Namely each single 
decimal digit occurs with frequency $1/10$, every pair of digits occurs with 
frequency $1/100$,... every block of $k$ digits occurs with frequency 
$1/10^k$,... Such a number is said to be ``normal to base $10$''. Borel asked 
in the 50's whether $\sqrt{2}$ is normal to base $10$: the answer is still 
unknown (also when replacing base $10$ with any integer base). Another hard 
to believe statement is that, although almost all real numbers are normal to 
any integer base, it is very difficult to provide an example of such a number.

\bigskip

Then note the following result and think about how to prove it: a real number 
$\xi$ is normal to base $b$ if and only if the sequence $(\xi b^n)_{n \geq 0}$ 
is uniformly distributed modulo $1$. Here a sequence $(x_n)_{n \geq 0}$ is said 
to be uniformly distributed modulo $1$ if the proportion of integers $n$
between $1$ and $N$ for which $\{x_n\}$ (the fractional part of $x_n$) belongs 
to any subinterval of $[0,1]$ is asymptotically equal to the length of this 
interval, i.e., if
$$
\forall u,v, \ \ 0 \leq u < v \leq 1, \ \ \lim_{n \to \infty} 
\frac{1}{N} \sharp \{n, \ 1 \leq n \leq N, \ u \leq \{x_n\} \leq v\} = v-u.
$$
This shows one of the reasons for studying sequences that are uniformly 
distributed modulo $1$. Unfortunately, not very much is known in this field 
for some specific sequences either. E.g., on one hand the sequence 
$(\xi \alpha^n)_{n \geq 0}$ is uniformly distributed modulo $1$ for all 
$\alpha > 1$ and almost all real numbers $\xi$, as well as for almost all 
$\alpha > 1$ and all nonzero real numbers $\xi$; on the other hand, 
frustratingly enough, one does not know whether the sequence 
$((\frac{3}{2})^n)_{n \geq 0}$ is uniformly distributed modulo $1$ (or even 
whether its limit points form a dense subset of the interval $[0,1)$) and the 
answer is unknown if $\frac{3}{2}$ is replaced by any rational number 
$\frac{p}{q} > 1$ such that $\frac{p}{q}$ is not an integer. Of course it is 
not a surprise that Pisot numbers enter this picture (recall that a Pisot
number is a real algebraic number $> 1$ such that all its conjugates lie 
inside the unit disk, hence the sequence of powers of a Pisot number tends to 
$0$ modulo $1$). May be less well-known is a link between the distribution of 
$((\frac{3}{2})^n)_{n \geq 0}$ modulo $1$ and the Waring problem. Namely let 
$g(n)$ be the smallest integer such that every positive integer can be expressed
as the sum of at most $g(n)$ $n$th powers ($n \geq 2$). Then we have 
$g(n) \geq 2^n + \lfloor (3/2)^n \rfloor - 2$. It can be proven that equality 
holds for $n$ if $\| (3/2)^n \| \geq (3/4)^{n-1}$ (where $\|x\|$ is the distance
between $x$ and its nearest integer neighbor). Note that Mahler proved that this
last inequality is true for $n$ large enough.

\noindent
Another both fascinating and frustrating question in distribution modulo $1$ is
the existence of the so-called Mahler $Z$-numbers: a positive integer $\xi$ is 
called a $Z$-number if for any $n \geq 0$ one has $0 \leq  \{\xi (3/2)^n \} < 
(1/2)$. A classical conjecture is that $Z$-numbers do not exist. Among the 
results in this direction, it can be proven that, if $b$ is an integer, then 
for any irrational number $\xi$ the numbers $\{\xi b^n \}$ cannot all lie in an 
interval of length $< 1/b$. Curiously enough equality is attained for ``Sturmian
numbers''. Recall that a number is Sturmian if its binary expansion is the 
coding of the trajectory with irrational initial slope of a ball bouncing on a 
square billiard, where bounces on vertical sides are coded by $0$ and bounces 
on horizontal sides are coded by $1$.

\bigskip

Changing the integer $b$ in some other integer base $b'$ in the study of 
$(\xi b^n)_{n \geq 0}$ involves comparing the expansion of a real number in two 
different bases. Again not much is known in the ``interesting case'' where the 
two bases $b$ and $b'$ are multiplicatively independent (i.e., 
$\frac{\log b'}{\log b}$ is irrational). We cite the Cassels-Schmidt result: 
{\it Let $r$ and $s$ be two multiplicatively independent integers $\geq 2$. 
Then the set of real numbers that are normal in base $r$ and not normal in 
base $s$ is uncountable}. Actually they prove more: the result still holds if 
one looks at the set of real numbers that are normal in base $r$, but not even 
simply normal in base $s$ (a number is said simply normal in base $s$ if each 
digit in $\{0, 1, \ldots, s-1\}$ occurs in it with frequency $1/s$).

\bigskip

We mention in passing that the above studies and questions are not only 
linked to transcendence of real numbers, but also to ``good approximations'' of 
reals by rationals. Furthermore, they are linked to the study of real numbers 
with missing digits in a given base, to questions about non-integer numeration 
bases (introduced by R\'enyi), etc. 

\bigskip

Before discussing continued fractions, we would like to point out two subjects 
that are addressed or alluded to in Bugeaud's book. The first one concerns real 
numbers $\xi$ whose base $b$ expansion is $d$-automatic. Recall that a sequence 
$(a_n)_{n \geq 0}$ is called $d$-automatic ($d$ integer $> 1$) if the set of 
subsequences $\{(a_{d^kn+j})_{n \geq 0}, \ k \geq 0, \ j \in [0, d^k-1]\}$ (the ``$d$-kernel'' of the sequence) is finite. Such sequences take necessarily only 
finitely many values; classical examples of nonperiodic $2$-automatic sequences 
are the Thue-Morse and the Shapiro-Rudin sequence. The Thue-Morse sequence can 
be defined by $a_n = s_2(n) \bmod 2$, where $s_2(n)$ stands for the sum of 
binary digits of the integer $n$. The Shapiro-Rudin sequence $(b_n)_{n \geq 0}$ 
can be defined by: $b_n$ is the number, reduced modulo $2$, of (possibly 
overlapping) occurrences of the block $11$ in the binary expansion of the 
integer $n$. Automatic sequences (the concept comes from theoretical computer 
science and can be given a more algorithmic definition) can be seen intuitively 
as sequences that can be computed by some sort of ``easy'' algorithm. One 
question about real numbers such that their base $b$ expansion is a 
$d$-automatic sequence is to prove that they cannot be algebraic irrational.
In other words, they are necessarily rational --the ``trivial'' case where the 
expansion is ultimately periodic-- or transcendental. Another formulation is 
that the digits of a number like $\sqrt{2}$ cannot be computed by a ``too simple
algorithm'' like a $d$-automaton). Mahler proved long ago that the Thue-Morse 
number is indeed transcendental; then Loxton and van der Poorten obtained 
partial results. The final result is due to Adamczewski and Bugeaud who used a 
clever mixture of Schmidt's subspace theorem and of combinatorial properties of 
``stammering sequences''. Note that the statement that automatic real numbers 
are necessarily rational or transcendental can be seen as a precise formulation 
of a particular case of a ``metatheorem'' stating that ``algebraic irrationals 
cannot have too simple base $b$ expansions''.  

\smallskip

The second subject we would like to cite in this section is the base $b$ 
expansion of real numbers for $b > 1$ not necessarily being integer. The author 
cites the fundamental result of Parry which gives a combinatorial 
characterization of the (closure of the) set of all possible $b$-expansions. 
Getting further in this direction would have lead the author to add another 
large amount of material. In particular, without entering details, we just want 
to recall that Erd\H{o}s, Jo\'o, and Komornik introduced, for a base $b \in 
(1,2)$, the set of real numbers between $0$ and $1$ that admit a unique 
expansion as $\sum a_n b^{-n}$ with $a_n \in \{0,1\}$. The $b$'s such that $1$ 
admits a unique base $b$ such expansion were later called ``univoque''. The 
reader will see the flavor of these numbers and the links with combinatorics of 
``words'' (i.e., finite sequences on a finite alphabet) if we add that univoque 
numbers $b$ can be characterized by the expansion of $1$: $b$ is univoque if and
only if the (unique) sequence $(a_n)_{n \geq 0}$ such that $1 = \sum a_n b^{-n}$
satisfies $(1-a_n)_{n \geq 0} < (a_{n+k})_{n \geq 0} < (a_n)_{n \geq 0}$ for all
$k > 0$, where $<$ is the lexicographical order on binary sequences.

\bigskip

Expanding a real number in some (possibly non-integer) base can be replaced by 
expanding a positive real number into a continued fraction. It can be proved for
example that almost all reals in $[0, 1]$ have a ``normal'' continued fraction 
expansion (for a reasonable definition of normality). But nothing is known about
the continued fractions of algebraic numbers of degree $\geq 3$: it is not 
known for example whether there exist algebraic numbers of degree $\geq 3$ 
having bounded partial quotients (actually it is not known either whether there 
exist algebraic numbers of degree $\geq 3$ having unbounded partial quotients;
these questions were asked by Khintchine). As above, one of the questions that 
can be asked is whether continued fractions having a $d$-automatic sequence of 
partial quotients must be either quadratic --this is the ``trivial case'' where 
the sequence of partial quotients is ultimately periodic-- or transcendental.

\bigskip

As it would be impossible to cite here all the results in the book, nor the nice
open questions given at the end, we recommend that the readers just open the 
book: they will, without noticing, jump from chapter to chapter and through 
(part of) the bibliography (of 751 items), eager to learn more about all these 
simple-to-state but hard-to-solve-or-still-open questions. We however cannot 
resist to give (only) five jewels: three are about the best results to date 
related to the questions of Borel and of Khintchine, two are open questions.

\medskip

1. [Adamczewski and Bugeaud] Let $\xi$ be a positive real number. Let 
$p_{\xi,b}(n)$ be the number of distinct blocks of digits of length $n$ in 
the base $b$ expansion of $\xi$. If $\xi$ is algebraic irrational, then
$$
\lim_{n \to \infty} \frac{p_{\xi,b}(n)}{n} = + \infty.
$$
A consequence of this result is that automatic real numbers are rational or 
transcendental.

\medskip

2. [Bugeaud and Evertse] Let $\xi$ be a positive real number. Let 
${\mathcal D}_{\xi,b}(n)$ be the number of digit changes up to the $n$th digit 
in the base $b$ expansion of $\xi$, i.e., if $\xi = \sum a_k b^{-k}$ with 
$a_k \in [0,b)$, then 
${\mathcal D}_{\xi,b}(n) = \sharp\{k, \ 1 \leq k \leq n, \ a_k \neq a_{k+1}\}$.
If $\xi$ is algebraic irrational, then there exists an effectively computable 
constant $n_0(\xi,b)$ depending only of $\xi$ and $b$ such that, for any integer
$n \geq n_0(\xi,b)$,
$$
{\mathcal D}_{\xi,b}(n) \geq (\log n)^{5/4}.
$$

\medskip

3. [Bugeaud] Automatic continued fractions are transcendental or quadratic. 
(An automatic continued fraction is a continued fraction whose sequence of 
partial quotients is $d$-automatic for some integer $d \geq 2$.)

\medskip

4. Open question [Furstenberg]. Suppose that the integers $p$ and $q$ are 
multiplicatively independent (i.e., they are not powers of the same integer). 
Then, in the expansion of $p^n$ to base $pq$, every digit and every combination 
of digits will occur, as soon as $n$ is sufficiently large.

\medskip

5. One of the open questions was asked by Mend\`es France (he attributes it to 
Mahler). Let $(c_n)_{n \geq 1}$ be a sequence of $0$'s and $1$'s. If the real 
numbers $\sum c_n 2^{-n}$ and $\sum c_n 3^{-n}$ are both algebraic, then they 
are both rational (i.e., the sequence $(c_n)_{n \geq 1}$ is ultimately 
periodic). This would be proved if another stronger question (Mahler) had an 
affirmative answer: is it true that the middle third Cantor set does not contain
any algebraic irrational number? Of course this last question would be answered 
affirmatively if, as Borel asked, not only $\sqrt{2}$ but also all algebraic 
irrational numbers would be proven to be normal in any base.

\end{document}